\documentclass[12pt]{amsart}
\usepackage{graphicx,verbatim,amssymb,amscd}
\usepackage[active]{srcltx}
\newtheorem{thm}{Theorem}

\theoremstyle{definition}

\theoremstyle{remark}

\numberwithin{equation}{section}

\def\CP{\mathbb{C}P}

\def\0{\varnothing}

\begin{document}

\title[Cut and semi-conjugate]
{Cut and semi-conjugate}
\author{Vladlen Timorin}

\address[Vladlen Timorin]{Faculty of Mathematics and Laboratory of Algebraic Geometry\\
National Research University Higher School of Economics\\
7 Vavilova St 117312 Moscow, Russia\\
Phone: +7 917 541 88 12}

\address[Vladlen~Timorin]
{Independent University of Moscow\\
Bolshoy Vlasyevskiy Pereulok 11, 119002 Moscow, Russia}

\email{vtimorin@hse.ru}

\thanks{
Partially supported by
the Deligne fellowship, the Simons-IUM fellowship, RFBR grants 10-01-00739-a,
11-01-00654-a, MESRF grant MK-2790.2011.1, and
AG Laboratory NRU-HSE, MESRF grant ag. 11 11.G34.31.0023
}


\begin{abstract}
We define a very general class of rational functions $f\,:\CP^1\to\CP^1$ 
such that for every function $f$ of this class, there exists a 
countable family of smooth curves $\gamma_i$ and a critically finite
function $R$ such that the dynamical system obtained from $f$ by
cutting along the curves $\gamma_i$ is topologically semi-conjugate to $R$.
\bigskip
\end{abstract}


\maketitle

We will consider topological dynamical systems on the 2-dimensional sphere $S^2$
given by continuous maps $f:S^2\to S^2$.
If such a map $f$ is fixed, it can be thought of as a
geometric structure on the sphere.
Namely, we can imagine that the sphere is equipped with arrows connecting
$x$ with $f(x)$ for all points $x\in S^2$.
Of course, the shape of the arrows does not matter, the only essential
information being which pairs of points are connected by arrows.

Homeomorphisms of spheres with arrows (i.e. topological conjugacies) 
to not change topological dynamics.
Thus, to change it, one needs some sort of topological surgery including
discontinuous operations such as cuts.
The simplest way to cut the sphere is to make a cut along some simple curve.
However, if the sphere comes with arrows depicting a dynamical system $f$,
then a cut along a curve creates some problems.
Namely, if we cut through the tip of an arrow, then the arrow splits.
This is problematic because a continuous map cannot create two different
arrows beginning at the same point.
To resolve this problem, one needs to make additional cuts, namely,
through the beginnings of all the arrows, whose ends are in the cut.
In other words, if we cut the sphere along a curve, then we also need
to cut it along the pre-image of this curve.
This creates further problems as we cut through the tips of some arrows again.
This, if we cut along a simple curve $Z\subset S^2$, then we also
need to cut along $f^{-1}(Z)$, along $f^{-2}(Z)$, and so on.
This leads to countably many cuts.
As we make all these cuts, we obtain some topological space $X$ called
a {\em sphere with cuts}.

The sphere with cuts $X$ carries well-defined arrows.
In other terms, we have a well-defined continuous map $F:X\to X$.
Note that the topological dynamical system $F:X\to X$ is uniquely
determined by the map $f$ and the first curve $Z$, along which we cut
(we call this curve the {\em initial cut}).
We will also consider a more general set-up, in which $Z$ is a union
of several simple curves.
We will refer to the map $F$ as the {\em map obtained from $f$ by cutting}.

More precisely, the set $X$ and the map $F$ are defined as follows.
Let $U_n$ be the complement to $Z\cup f(Z)\cup\dots\cup f^{\circ n}(Z)$ in $S^2$.
Define the compact space $X_n$ as the Caratheodory compactification of $U_n$
(i.e. the union of $U_n$ and the set of all prime ends of $U_n$ equipped with
a suitable topology).
The inclusion $U_{n+1}\to U_n$ gives rise to a continuous map $\iota_n:X_{n+1}\to X_n$. 
We define the topological space $X$ as the inverse limit of the system of spaces $X_n$ and maps $\iota_n$.
The maps $F_n:X_{n+1}\to X_n$ corresponding to $f:U_{n+1}\to U_n$ define 
the continuous map $F:X\to X$.

After countably many cuts, the orbits of almost all points remain unchanged.
Thus we can think that the dynamics of the map $F:X\to X$ is not too
much different from the dynamics of the map $f$.
On the other hand, as a topological space, $X$ can be very different from
the sphere, e.g. it can have uncountably many connected components.
In a very general situation described below, there is a semi-conjugacy
between the topological dynamical system $F:X\to X$ and some {\em hyperbolic
critically finite} (i.e. very good!) function $R$.
Recall that a function $R$ is hyperbolic and critically finite if each
of its critical points gets eventually mapped to a periodic critical point.
Thus we can study the dynamics of $F:X\to X$
(hence also the dynamics of $f:\CP^1\to\CP^1$) using a semi-conjugacy with
a hyperbolic rational function.
In fact, the good rational function $R$ can be usually chosen in infinitely
many ways.

To state a precise theorem, we will need a relation on the
set of all rational functions.
It resembles combinatorial equivalence but is in fact much weaker.
Let $f$ and $R$ be rational functions, and suppose that $R$ has a
finite post-critical set $P_R$.
We say that the function $R$ {\em can represent} the function $f$, if 
$R$ is homotopic to a map topologically conjugate to $f$ through branched
coverings so that the homotopy preserves the finite set $R(P_R)$.  
The difference with the usual definition of combinatorial equivalence is
only that the set $P_R$ is replaced with the set $R(P_R)$.
But, in contrast to combinatorial equivalence, there are usually many
(even infinitely many) critically finite rational functions $R$
that can represent the same function $f$.
However, if $R$ can represent $f$, then these two functions must have
the same structure of super-attracting orbits.
On the other hand, the function $f$ can have many critical points,
whose orbits are infinite and show as complicated behavior as they please.
We require that $R$ be hyperbolic.
In particular, it must have at least one super-attracting cycle.
This means that $f$ must also have at least one super-attracting cycle; 
moreover, the structure of super-attracting cycles for $f$ must be the same as for $R$.
This is almost all we need from $f$.

We now state the main result.

\begin{thm}
  \label{t:semiconj}
  Let $f:\CP^1\to\CP^1$ be a rational function, and $R:\CP^1\to\CP^1$
  a critically finite hyperbolic function that can represent $f$.
  Then there exist an initial cut $Z$ (which is a finite union of
  smooth simple curves), the corresponding sphere with cuts $X$ and
  the map $F:X\to X$ obtained from $f$ by cutting, such that the
  topological dynamical system $F:X\to X$ is semi-conjugate to the dynamical system
  $R:\CP^1\to\CP^1$.
\end{thm}

Recall that the dynamical system $F:X\to X$ is determined by
the choice of the initial cut $Z$.
This choice can be made very explicit.
Moreover, the set $Z$ is defined only up to a suitable homotopy,
so that we can arrange that all curves that appear in $Z$
are analytic, or even real semi-algebraic (alternatively, we can
make them broken lines).

Below, we give some examples, where the choice of $Z$ is made more explicit.
To this end, we consider certain complex one-dimensional parameter spaces
of quadratic rational functions.

It is a general observation and philosophy that the dynamical behavior
of a rational function is determined by the the behavior of
its critical orbits.
A rational function of degree two has two critical points.
Thus, if we fix a simple dynamical behavior of one of these points,
then we are left with only one ``free'' critical point.
This makes the corresponding parameter space complex one-dimensional.
To fix a simple behavior of a critical point, we can e.g. require
that this point be periodic of period $k$.
The space of such rational functions is $Per_k(0)$ (in this notation,
$0$ stands for the multiplier of a $k$-periodic point).
More precisely, $Per_k(0)$ consists of M\"obius conjugacy classes of
rational functions $R$ of degree 2 with marked critical points $c_1$,
$c_2$ such that $R^{\circ k}(c_1)=c_1$.

The parameter spaces $Per_k(0)$ are one-dimensional slices of
the parameter space of all degree 2 rational functions (with marked critical points).
E.g., for $k=1$, we obtain the space of quadratic polynomials $z^2+c$.
Indeed, if a rational function has a fixed critical point, then,
mapping this point to infinity by a suitable M\"obius conjugacy,
we can make this function into a quadratic polynomial; the
term with $z$ can be killed by a parallel translation.

The set of hyperbolic functions in $Per_k(0)$ is an open set, whose
components are called {\em hyperbolic components}.
Recall that a hyperbolic function $R$ with periodic critical point $c_1$
(and the corresponding hyperbolic component in $Per_k(0)$) is said to be of
{\em type C} (C stands for ``capture''; not to be confused with the capture operation
of B. Wittner) if $c_2$ lies in the basin of the the cycle $c_1$, $\dots$, $R^{\circ k-1}(c_1)$, 
but not in the immediate basin.
Recall also that any hyperbolic component of type C in $Per_k(0)$ contains a unique
critically finite M\"obius conjugacy class of rational functions.
In $Per_k(0)$, {\em any hyperbolic critically finite function of type C can represent any
non-hyperbolic function}.

So let $f$ be any non-hyperbolic function and $R$ any critically finite type C
hyperbolic function, whose classes are in $Per_k(0)$.
{\em The initial cut $Z$ is a simple curve containing the non-periodic 
critical point of $f$.}
It is defined as the full pre-image of some simple curve $\beta$ under $f$,
so that the map $f|_Z$ is two-to-one except that the critical point $c_2$ 
is the only preimage of the critical value $v=f(c_2)$.  
The curve $\beta$ connects the non-periodic critical value $v$ of $f$ with a certain
pre-periodic point $w$ that gets eventually mapped to the periodic critical
point of $f$.
{\em The curve $\beta$ is only defined up to a homotopy relative to the forward
orbit of $w$}.
The following is a useful property that characterizes the homotopy class of $\beta$.
Let $\sigma_\beta$ be a self-homeomorphism of the Riemann sphere that is equal
to the identity outside a small neighborhood of $\beta$ and that takes $v$ to $w$.
Then $\sigma_\beta\circ f$ is a critically finite branched covering, whose
combinatorial equivalence class is well defined.
The curve $\beta$ should be chosen so that to make the covering $\sigma_\beta\circ f$
combinatorially equivalent to $R$.
There is a simple and rather explicit construction of M. Rees \cite{R}
that produces such a curve $\beta$.

Suppose that the class of $f$ lies on the boundary of the type C hyperbolic 
component containing the class of $R$. 
Then $\beta$ can be chosen as the closure of an internal ray 
in a Fatou component of $f$, whose boundary contains the critical value $v$. 
In this case, all cuts we need to make in the sphere are disjoint simple 
curves.  
To obtain the topological dynamics of $R$ from the topological dynamics of $f$,
we need to cut along all pullbacks of $Z$ and then reglue them differently
(this is similar to cutting the plane along the interval $[-1,1]$, opening up
the cut to form the circle $x^2+y^2=1$, so that the previously
identified points have equal $x$-coordinates, and gluing the circle 
so that the points with equal $y$-coordinates are identified).
This topological surgery is called {\em regluing}; it was introduced in \cite{T09}. 
Regluing is reversible: if $R$ can be obtained from $f$ by regluing, then
also $f$ can be obtained from $R$ by regluing. 
This provides topological models for non-hyperbolic rational functions,
corresponding to boundary points of type C hyperbolic components in $Per_k(0)$.

A proof of Theorem \ref{t:semiconj} is given in preprint \cite{T10}.
The proof is rather technical; it uses a version of Thurston's algorithm
\cite{DH93}.
Thurston's algorithm is applied to the composition of $F:X\to X$
with the projection $X\to\CP^1$ rather than to a branched covering of the sphere.

\end{document}